\documentclass[twoside]{article}

\catcode`\@=11
\long\def\@makefntext#1{
\protect\noindent \hbox to 3.2pt {\hskip-.9pt  
$^{{\eightrm\@thefnmark}}$\hfil}#1\hfill}		

\def\@makefnmark{\hbox to 0pt{$^{\@thefnmark}$\hss}}	
	
\def\ps@myheadings{\let\@mkboth\@gobbletwo
\def\@oddhead{\hfill\hbox{}\rightmark}   
\def\@oddfoot{}\def\@evenhead{\leftmark\hbox{}\hfill}\def\@evenfoot{}
\def\sectionmark##1{}\def\subsectionmark##1{}}

\oddsidemargin=\evensidemargin
\newcounter{sectionc}\newcounter{subsectionc}\newcounter{subsubsectionc}
\renewcommand{\section}[1] {\vspace{12pt}\addtocounter{sectionc}{1} 
\setcounter{subsectionc}{0}\setcounter{subsubsectionc}{0}\noindent 
	{\tenbf\thesectionc. #1}\par\vspace{5pt}}
\renewcommand{\subsection}[1] {\vspace{12pt}\addtocounter{subsectionc}{1} 
	\setcounter{subsubsectionc}{0}\noindent 
  {\bf\thesectionc.\thesubsectionc. {\kern1pt \bfit #1}}\par\vspace{5pt}}
\renewcommand{\subsubsection}[1] {\vspace{12pt}\addtocounter{subsubsectionc}{1}
	\noindent{\tenrm\thesectionc.\thesubsectionc.\thesubsubsectionc.
	{\kern1pt \tenit #1}}\par\vspace{5pt}}
\newcommand{\nonumsection}[1] {\vspace{12pt}\noindent{\tenbf #1}
	\par\vspace{5pt}}

\newcounter{appendixc}
\newcounter{subappendixc}[appendixc]
\newcounter{subsubappendixc}[subappendixc]
\renewcommand{\thesubappendixc}{\Alph{appendixc}.\arabic{subappendixc}}
\renewcommand{\thesubsubappendixc}
	{\Alph{appendixc}.\arabic{subappendixc}.\arabic{subsubappendixc}}

\renewcommand{\appendix}[1] {\vspace{12pt}
        \refstepcounter{appendixc}
        \setcounter{figure}{0}
        \setcounter{table}{0}
        \setcounter{lemma}{0}
        \setcounter{theorem}{0}
        \setcounter{corollary}{0}
        \setcounter{definition}{0}
        \setcounter{equation}{0}
        \renewcommand{\thefigure}{\Alph{appendixc}.\arabic{figure}}
        \renewcommand{\thetable}{\Alph{appendixc}.\arabic{table}}
        \renewcommand{\theappendixc}{\Alph{appendixc}}
        \renewcommand{\thelemma}{\Alph{appendixc}.\arabic{lemma}}
        \renewcommand{\thetheorem}{\Alph{appendixc}.\arabic{theorem}}
        \renewcommand{\thedefinition}{\Alph{appendixc}.\arabic{definition}}
        \renewcommand{\thecorollary}{\Alph{appendixc}.\arabic{corollary}}
        \renewcommand{\theequation}{\Alph{appendixc}.\arabic{equation}}
        \noindent{\tenbf Appendix#1}\par\vspace{5pt}}
\newcommand{\subappendix}[1] {\vspace{12pt}
        \refstepcounter{subappendixc}
        \noindent{\bf Appendix \thesubappendixc. {\kern1pt \bfit #1}}
	\par\vspace{5pt}}
\newcommand{\subsubappendix}[1] {\vspace{12pt}
        \refstepcounter{subsubappendixc}
        \noindent{\rm Appendix \thesubsubappendixc. {\kern1pt \tenit #1}}
	\par\vspace{5pt}}

\newcommand{\textlineskip}{\baselineskip=13pt}
\newcommand{\smalllineskip}{\baselineskip=10pt}

\newcommand{\copyrightheading}[1]
	{\vspace*{-2.5cm}\smalllineskip{\flushleft
	{\footnotesize #1}
	 }}

\def\abstracts#1#2#3{{
	\centering{\begin{minipage}{4.5in}\footnotesize\baselineskip=10pt
	\parindent=0pt #1\par 
	\parindent=15pt #2\par
	\parindent=15pt #3
	\end{minipage}}\par}} 

\def\keywords#1{{
	\centering{\begin{minipage}{4.5in}\footnotesize\baselineskip=10pt
	{\footnotesize\it Keywords}\/: #1
	 \end{minipage}}\par}}

\newcommand{\bibbf}{\ninebf}
\renewenvironment{thebibliography}[1]
	{\frenchspacing
	 \ninerm\baselineskip=11pt
	 \begin{list}{\arabic{enumi}.}
        {\usecounter{enumi}\setlength{\parsep}{0pt}     
	 \setlength{\leftmargin 12.7pt}{\rightmargin 0pt} 
         \setlength{\itemsep}{0pt} \settowidth
	{\labelwidth}{#1.}\sloppy}}{\end{list}}

\newcounter{itemlistc}
\newcounter{romanlistc}
\newcounter{alphlistc}
\newcounter{arabiclistc}

\newenvironment{alphlist}
	{\setcounter{alphlistc}{0}
	 \begin{list}{$($\alph{alphlistc}$)$}
	{\usecounter{alphlistc}
	 \setlength{\parsep}{0pt}
	 \setlength{\itemsep}{0pt}}}{\end{list}}

\newcommand{\fcaption}[1]{
        \refstepcounter{figure}
        \setbox\@tempboxa = \hbox{\footnotesize Fig.~\thefigure. #1}
        \ifdim \wd\@tempboxa > 5in
           {\begin{center}
        \parbox{5in}{\footnotesize\smalllineskip Fig.~\thefigure. #1}
            \end{center}}
        \else
             {\begin{center}
             {\footnotesize Fig.~\thefigure. #1}
              \end{center}}
        \fi}

\newcommand{\tcaption}[1]{
        \refstepcounter{table}
        \setbox\@tempboxa = \hbox{\footnotesize Table~\thetable. #1}
        \ifdim \wd\@tempboxa > 5in
           {\begin{center}
        \parbox{5in}{\footnotesize\smalllineskip Table~\thetable. #1}
            \end{center}}
        \else
             {\begin{center}
             {\footnotesize Table~\thetable. #1}
              \end{center}}
        \fi}

\def\pmb#1{\setbox0=\hbox{#1}
	\kern-.025em\copy0\kern-\wd0
	\kern.05em\copy0\kern-\wd0
	\kern-.025em\raise.0433em\box0}

\def\fnt#1#2{\footnotetext{\kern-.3em
        {$^{\mbox{\scriptsize #1}}$}{#2}}}

\font\tenrm=cmr10
\font\tenit=cmti10 
\font\tenbf=cmbx10
\font\bfit=cmbxti10 at 10pt
\font\ninerm=cmr9

\font\ninebf=cmbx9
\font\eightrm=cmr8

\def\qed{\hbox{${\vcenter{\vbox{		 
   \hrule height 0.4pt\hbox{\vrule width 0.4pt height 6pt
   \kern5pt\vrule width 0.4pt}\hrule height 0.4pt}}}$}}

\def\theequation{\thesectionc.\arabic{equation}}  

    	{\setcounter{itemlistc}{0}		
	 \begin{list}{$\bullet$}		
	{\usecounter{itemlistc}			
	 \leftmargin10pt	       
	 \setlength{\parsep}{0pt}
	 \setlength{\itemsep}{0pt}     
	}}{\end{list}}

	{\setcounter{romanlistc}{0}		
	 \begin{list}{$($\roman{romanlistc}$)$}	
	{\usecounter{romanlistc}		
	 \leftmargin18pt 
	 \setlength{\parsep}{0pt}
	 \setlength{\itemsep}{0pt}	
	 \settowidth{\labelwidth}{#1}                          
	}}{\end{list}}

	{\setcounter{enumii}{0}			
	 \begin{list}{$($\alph{enumii}$)$}	
	{\usecounter{enumii}			
	 \leftmargin18pt		
	 \setlength{\parsep}{0pt}
	 \setlength{\itemsep}{0pt}	
	 \settowidth{\labelwidth}{#1}                          
	}}{\end{list}}

\def\bsc{{\sc a\kern-6.4pt\sc a\kern-6.4pt\sc a}}  	
\def\bflatex{\bf L\kern-.30em\raise.3ex\hbox{\bsc}\kern-.14em 
T\kern-.1667em\lower.7ex\hbox{E}\kern-.125em X} 

\usepackage{amsfonts}

\newcommand{\C}{{\mathbb{C}}}
\newcommand{\R}{{\mathbb{R}}}
\newcommand{\Q}{{\mathbb{Q}}}
\newcommand{\N}{{\mathbb{N}}}
\newcommand{\D}{{\mathbb{D}}}

\def\cd{\widehat{\C}}
\def\BI{I}
\def\id{\mbox{id}}
\def\dist{\mbox{dist}}

\begin{document}
\setlength{\textheight}{7.7truein}    

\markboth{\protect{\footnotesize\it M. Wendt}}{\protect{\footnotesize\it The entropy of entire transcendental functions}}

\normalsize\textlineskip
\thispagestyle{empty}
\setcounter{page}{1}

\copyrightheading{}

\vspace*{1in}         

\centerline{\bf THE ENTROPY OF ENTIRE TRANSCENDENTAL FUNCTIONS}
\vspace*{0.37truein}
\centerline{\footnotesize MARKUS WENDT\footnote{Mathematisches Seminar, Christian-Albrechts-Universit\"at zu Kiel,Ludewig-Meyn-Str. 4, D-24098 Kiel, Germany }}


\vspace*{0.21truein}
\abstracts{We use Bowen's definition of topological entropy 
  and Ahlfors five islands theorem, as well as the theory of polynomial-like mappings, 
  to show that the topological entropy of any entire transcendental
  function is infinity. In addition the entropy is 
  concentrated  on the Julia set for each meromorphic function which has 
  no wandering domains.}{}{}

\vspace*{5pt}

\vspace*{4pt}
\baselineskip=13pt              
\normalsize              	
\section{Introduction}          
\noindent
 If $X$ is a topological space, $A\subseteq X$ and $f:X\to X$ continuous,
 we write $h(f,A)$ for the topological entropy of $f$ on $A$
 defined  by Bowen \cite{Bowen:73}. 
 Then, the topological entropy of $f$ is given by $h(f) := h(f,X)$.
 Bowen's definition also works for arbitrary self-mappings,
 only a few statements have to be modified. We start with 
 a short introduction to Bowen's definition of topological entropy.
  
 \vspace*{12pt}
 \noindent
 {\bf Notation~}  
  Let $X$ be a topological space, $C\subseteq X$  and let $\xi$ be a cover 
  of $X$.
  We write $C\prec \xi$ if there exists a $E\in \xi$ with $C\subseteq E$. 
  Let $f:X\to X$ be a map.
  We set $f^0 := \id$, $f^n := f \circ f^{n-1}$, $n\in \N$, and $f^{-n} := (f^n)^{-1}$. 
  Let
  \[
     N_{f,\xi}(C) := \left\{ \begin{array}{ll}
                       \sup\bigl\{ k \in \N_0 : f^\ell(C) \prec \xi \hbox{ for $\ell =0, \ldots, k$}\bigr\}, & \hbox{ if $C\prec \xi$} \\[0.05cm]
                        0, & \hbox{ otherwise}
                     \end{array} \right. 
  \]
   and $D_{f,\xi}(C) :=  \exp\bigl( -N_{f,\xi}(C) \bigr)$ with  $\exp(-\infty) := 0$. 

 \vspace*{12pt}
 \noindent
 {\bf Definition 1.1~}  
   Let $X$ be a topological space, $A\subseteq X$,  $f:X\to X$
    and let $\xi$ be a cover of $X$. We set for all $s, \delta \in \, ]0,\infty[$
   \[
      h_{\xi, \delta}^s (f,A) := \inf \Bigg\{ \sum\limits_{j=1}^\infty D_{f,\xi}(C_j) ^s : A \subseteq \bigcup\limits_{j \in \N} C_j, \;D_{f,\xi}(C_j)  \leq \delta \hbox{ for all $j\in \N$}\Bigg\}
   \]
   and let $h_\xi^s(f,A) :=  \sup\limits_{\delta > 0} h_{\xi, \delta}^s(f,A) = \lim\limits_{\delta \to 0} h_{\xi, \delta}^s(f,A)$,
   as well as
   \[
     h_\xi (f,A) :=  \inf\bigl\{ s \in \, ]0,\infty[:\; h_\xi^s(f,A) = 0 \bigr\}  .
   \]
   The {\bf topological entropy of $f$ on $A$}  is defined by
    \[
     h(f,A) := \sup \bigl\{ h_\xi(f,A): \xi \hbox{ finite open cover of $X$} \bigr\},
    \]
     and  $h(f) :=  h(f,X)$ is called {\bfseries topological entropy} of $f$.

 \vspace*{12pt}
  \noindent
  It is well known \cite[Prop. 1]{Bowen:73}, \cite{Misiurewicz:2003} that this definition coincide with the usual 
  definiti\-ons given in \cite{Walters:82},
  if $X$ is a compact metric space and if $f$ is continuous.
  Almost all statements of \cite[Prop. 2]{Bowen:73} can be proven in this 
  general setting, the proofs are elementary \cite{Wendt:2002}.

 \vspace*{12pt}
 \noindent
 {\bf Proposition~1.2} (Properties of the topological entropy) {\it Let $X,X_1, X_2$ 
   be topo\-logical spaces and let $f:X\to X$,
           $f_1: X_1 \to X_1$, $f_2:X_2 \to X_2$.
   \begin{alphlist}
     \item If $f_1$ and $f_2$ are topologically conjugate,
                 i.e. there exists a homeomorphism $\pi:X_1 \to X_2$ with $\pi \circ f_1 = f_2 \circ \pi$,
                 then $h(f_1, A) = h(f_2, \pi(A))$ for all $A\subseteq X_1$.

    \item It holds $h\bigl(f,\bigcup\limits_{n\in \N} A_n\bigr) = \sup\limits_{n\in \N} h(f,A_n)$ for each sequence $(A_n)$ of subsets of $X$.		 
    \item For all $n\in \N_0$, $A\subseteq X$, we have $h(f^n,A) \leq n \cdot h(f,A)$ with equality if $f$ is continuous.
    
    \item  For each closed subset $A$ of $X$ with $f(A) \subseteq A$, 
           we have $h(f,A) = h(f\vert_A)$, where $f\vert_A$  denotes the restriction of $f$ to $A$.
    
    \item  For all $A \subseteq X$, we have $h(f,A) \leq h(f,f(A))$ 
           with equality if $f$ is continuous.
   \end{alphlist}}

 \vspace*{2pt}
 \noindent
 If we interpret strict positive entropy as a kind of chaos, we expect
 the following lemma, where no continuity is necessary.
 
 \vspace*{12pt}
 \noindent
 {\bf Lemma~1.3} {\it Let $X$ be a topological space, $A\subseteq X$, $x_0 \in X$
  and let $f:X\to X$ with $f^n\vert_A \to x_0$ (uniformly).
 Then $h(f,A) = 0$.}
 
 \vspace*{12pt}
  \noindent    
 {\bf Proof.} Let $\xi$ be a finite open cover of $X$ and let $\xi_n := \bigvee_{k=0}^n f^{-k} \xi$, $n\in \N$,
   being the refinement of $\xi$ with respect to $\id, f^{-1}, \ldots, f^{-n}$,
   i.e. $\xi_n = \bigl\{ \bigcap_{k=0}^n f^{-k}(E_k): E_1, \ldots, E_n \in \xi \bigr\}$.   
   Choose same $E^* \in \xi$ with $x_0 \in E^*$. The uniformly convergence implies
   the existence of same $n \in \N$ with
   \[
    f^m(A) \subseteq E^* \quad \hbox{for all $m\in \N_{\geq n}$.}
   \]
   For each $C \in A \cap \xi_n := \{ A \cap E: E \in \xi_n\}$ we have
   $f^\ell(C) \prec \xi$ for  all $\ell \in \N_0$, 
   i.e. $D_{f,\xi}(C) = 0$. It follows 
   $h_{\xi,\delta}^s(f,A) \leq \sum_{ C \in A \cap \xi_n} D_{f,\xi}(C)^s = 0$   
   for all $s,\delta \in \, ]0,\infty[$, thus $h_\xi(f,A) = 0$,  hence $h(f,A) = 0$. \, \qed 

 \vspace*{12pt}
 \noindent
  The next Proposition shows that the entropy of a map
  doesn't change if we add finitely many points to its domain.

 \vspace*{12pt}
 \clearpage
 \noindent
 {\bf Proposition~1.4} {\it Let $X'$ be a metric space, $X\subseteq X'$
 and let $f:X\to X$. Let $F:X'\to X'$ be a map with $F\vert_X = f$. Then
 \[
    h(F,A) \leq h(f,A) \quad \hbox{for all $A\subseteq X$}
 \]
 with equality if  $X'\setminus X$ is finite.}
 
 \vspace*{12pt}
  \noindent    
 {\bf Proof.}  
   Let $A\subseteq X$. "$\leq$": Let $\xi'$
   be a finite open cover of $X'$. Consider $\xi := \{ E \cap X: E \in \xi' \}$.
   Then $\xi$ is a finite cover of $X$, consisting of $X$-open sets. Let $s\geq 0$ with $h_\xi^s(f,A) = 0$.
   Let $\epsilon,\delta > 0$. There exists a cover $(C_j)_{n\in \N}$ of $A$ in $X$
   with $D_{f,\xi}(C_j) \leq\delta$ for all $j\in \N$ and $\sum_{j=1}^ \infty D_{f,\xi}(C_j)^s < \epsilon$.
   Thus $D_{F,\xi'}(C_j) \leq D_{f,\xi}(C_j) \leq \delta$ and $h_{\xi',\delta}^s(F,A) \leq \sum_{j=1}^\infty D_{F,\xi'}(C_j)^s \leq \sum_{j=1}^\infty D_{f,\xi}(C_j)^s < \epsilon$,
   hence $h_{\xi'}^s(F,A) = 0$, i.e. $h(F,A) \leq h(f,A)$.
   "$\geq$": Let $X'\setminus X = \{x_1, \ldots, x_k\}$.
   Choose open neighbourhoods $U_1, \ldots, U_k$ of 
   $x_1, \ldots, x_k$ with $\overline{U_i} \cap \overline{U_j} = \emptyset$
   for $i\not=j$.
   Let $\xi=\{ X \cap E_1, \ldots, X \cap E_m\}$ be a finite open cover of 
   $X$ consisting of $X$-open sets. Set
   \[
     \xi' := \{E_1, \ldots, E_m\} \cup \{ U_1, \ldots, U_k \}.
   \]
   Thus $\xi'$ is a finite open cover of $X'$. We
   show $h_{\xi'}(F,A) \geq h_\xi(f,A)$. 
   Let $s\geq 0$ with $h_{\xi'}^s(F,A) = 0$. It is 
   sufficient  to show $h_\xi^s(f,A) = 0$. Let $\epsilon, \delta > 0$.
   Then there exists a sequence $(C_j')_{j\in \N}$ with $C_j' \subseteq X'$,
   $\bigcup_{j\in \N} C_j' = A$ and $D_{\textstyle _{F,\xi'}}( C_j') \leq \delta$ for all $j\in \N$,
   as well as $\sum_{j=1}^ \infty D_{\textstyle _{F, \xi'}}(C_j') < \frac{\epsilon}{m}$.
   For all $j\in \N$ set $C_j := C_j' \cap X$ and
   \[
      L^j := \bigl\{ \ell \in \N_0:  F^\ell(C_j') \prec \xi', \; f^\ell(C_j) \not\prec \xi \bigr\}.
   \]
   Let $j\in \N$ with $L^j \not= \emptyset$ and let $\ell \in L^j$. 
   It follows $F^\ell(C_j') \subseteq U_{i_\ell}$
   for same $i_\ell \in \{1, \ldots, k\}$,
   thus $f^\ell(C_j) \subseteq U_{i_\ell} \setminus \{x_{i_\ell}\}$.
   Set
   \[
     H_{j,n}^\ell := f^{-\ell}\bigl( U_{i_\ell} \cap E_n \bigr) \cap C_j, \qquad H_{j,n} := \bigcap_{\ell \in L^j} H_{j,n}^\ell \quad \hbox{for $n=1, \ldots, m$.}
   \]
   It follows $C_j = H_{j,1} \cup \cdots \cup H_{j,m}$ and
   if $\ell \in \N_0$ with $F^\ell(C_j') \prec \xi'$, 
   we have $f^\ell(H_{j,n}) \prec \xi$ for $n=1, \ldots, m$.
   Thus $D_{f,\xi}(H_{j,n}) \leq D_{\textstyle_{F, \xi'}}(C_j') \leq \delta$,
   and
   \begin{eqnarray*}
      h_{\xi,\delta}^s(f,A) &\leq& \sum_{j=1,\atop L^j=\emptyset}^\infty D_{f,\xi}(C_j)^s  + \sum_{j=1,\atop L^j\not=\emptyset}^\infty \sum_{n=1}^m D_{f,\xi}(H_{j,n})^s  \\
       &\leq&      \sum_{j=1,\atop L^j=\emptyset}^\infty D_{F,\xi'}(C_j')^s  + \sum_{j=1,\atop L^j\not=\emptyset}^\infty m \cdot D_{F,\xi'}(C_j')^s   \leq m \cdot \sum_{j=1}^\infty D_{F, \xi'}(C_j')^s <  \epsilon,
   \end{eqnarray*}
   hence $h_{\xi}^s(f,A) = 0$. Thus $h(F,A) \geq h(f,A)$. \, \qed

\section{Zero entropy and the Fatou set}
 \noindent
 We write $\cd = \C \cup\{\infty\}$ for the Riemannian sphere and let
 \[
    \BI(X) := \{ f:X\to X\vert  \hbox{ meromorphic and \textit{not} bijective} \}
 \]
 for $X=\cd,\C,\C\setminus\{0\}$, and $\BI = \BI(\cd) \cup \BI(\C) \cup \BI(\C\setminus \{0\})$.
 For every map $f\in \BI$ we write $J(f)$ and $F(f)$ for the Julia respectively Fatou set
 of $f$ \cite{Bergweiler:93}. See  \cite{Hua:98,Milnor:2000} for an introduction to complex dynamic,
 especially for the classification theorem of periodic components of
 the Fatou set.
  The main goal of this section is Theorem 2.4. We show, that the
  topological entropy of a map $f\in \BI$ is concentrated on the 
  Julia set, namely $h(f,F(f))=0$, if $f$ has no wandering domain.
  
\vspace*{12pt}
\noindent
{\bf Lemma~2.1.} {\it Let $f \in \BI$ and let $U$ be a B\"ottcher- Schr\"oder or Leau domain
 of $f$. Then $ h(f,U) = 0$.}
 
\vspace*{12pt}
\noindent
{\bf Proof.} Choose the corresponding fixpoint $z$ of $f$ with respect to $U$.
 Choose a sequence $(A_m)$ of compact sets  with 
 $U=\bigcup_{m\in \N} A_m$. Let $m\in \N$.  It follows $f^n\vert A_m \to z$ uniformly. 
 Lemma 1.3 implies $h(f,A_m) = 0$, and  by Proposition 1.2 we get
 \[
   h(f,U) = h(f, \bigcup_{m\in \N} A_m) = \sup_{m\in \N} h(f,A_m) = 0.     \qquad \qquad  \qed 
 \]
 
 \vspace*{2pt}
 \noindent
 {\bf Lemma~2.2.} {\it Let $f \in \BI$ and let $U$ be a Siegel disc or a Herman ring of $f$.
    Then $h(f,U) = 0$.}
    
 \vspace*{12pt}
 \noindent    
 {\bf Proof.}  Let $U$ be a Siegel disc. Thus there exists a bijective, meromorphic 
 map $\phi: \D \to U$, where $\D=\{z \in \C: \vert z \vert < 1\}$, and
 $\alpha \in \R\setminus \Q$ such that $\phi(\lambda \cdot z) = f( \phi(z))$
 for $\lambda := e^{2\pi i \alpha}$ and all $z\in \D$. Hence $g:\D \to \D, z\mapsto \lambda \cdot z$
 and $f\vert_{U}:U\to U$ are conjugated, thus $h(f\vert_U,U) = h(g,\D)$.
 Consider $A_n := \{ z \in \C: \vert z \vert \leq 1- \frac{1}{n} \}$, $n\in \N$.
 Using Proposition 1.2 and Proposition 1.4, we get 
 \[
    h(f,U) \leq   h(f\vert_U, U ) = h(g,\D) = h(g,\bigcup_{n\in \N} A_n) = \sup_{n\in \N} h(g,A_n) = \sup_{n\in \N} h(g\vert A_n) = 0,
 \]
  where the last equation follows from the fact, that the entropy of every rotation 
  is zero \cite{Walters:82}. The case of a Herman ring is similar.  \, \qed

 \vspace*{12pt}
 \noindent
 {\bf Lemma~2.3.} {\it Let $f \in \BI$ and let $U$ be a Baker domain of $f$.
    Then $h(f,U) = 0$.}
    
 \vspace*{12pt}
 \noindent    
 {\bf Proof.}  We have $f^n\vert_U \to z_0$
  where $z_0=0$ or $z_0=\infty$  doesn't lie in the domain of $f$. 
   Consider  $F:\cd \to \cd$, defined by 
  \[
    F(z) := \left\{ \begin{array}{ll}
    	      f(z), & \hbox{ if $f(z)$ is well defined} \\
              \infty, & \hbox{ if $f(z)$ is not well defined, and $z=\infty$} \\
	      0, & \hbox{ otherwise}
           \end{array}\right.
  \]
  Now we can apply Lemma 1.3 to $F$ in the same way as in the proof of
  Lemma 2.1. It follows $0 = h(F, U) = h( f, U)$ by Proposition 1.4. \, \qed

 \vspace*{12pt}
 \noindent
 {\bf Theorem~2.4.} {\it Let $f \in \BI$ without having 
       wandering domains. Then $h(f, F(f)) = 0$, especially $h(f, J(f))= h(f)$.}

 \vspace*{12pt}
  \noindent    
 {\bf Proof.} Let $A\subseteq F(f)$ compact. For every $x\in A$ there
 exists a pre-periodic component $U_x$ of $F(f)$ with $x\in U_x$. It follows
 $A\subseteq \bigcup_{j=1}^n U_{x_j}$ for same $n\in \N, $ $x_1, \ldots, x_n \in A$. 
 Let $j \in \{1, \ldots, n\}$ and $U=U_{x_j}$. There exists a periodic component
 $W$ of $F(f)$ with $f^m(U) \subseteq W$ for same $m\in \N$. Let $p\in \N$ be
 the period of $U$. The classification theorem implies that $W$ is a 
 B\"ottcher-, Schr\"oder-, Leau-, Baker domain, a Siegel disc or a Herman ring (for $f^p$). 
 Using some properties of the topological entropy (Proposition 1.2) 
 and Lemma 2.1-2.3, it follows
 \[
   0=h(f^p, W) = p \cdot h(f,W) \geq p \cdot h(f, f^m(U)) = p \cdot h(f,U) \geq 0,
 \]
 hence $h(f^p,U ) = 0$. It follows $h(f,A) = \max_{j=1}^ n h(f,U_{x_j}) = 0$. Using
 a sequence $(A_n)$ of compact sets with $\bigcup_{n\in \N} A_n = F(f)$, we have
 shown $h(f,A_n) = 0$ for every $n\in \N$, thus $h(f,F(f) ) = \sup_{n\in \N} h(f,A_n) = 0$. \, \qed

\section{Graph theory and the Ahlfors five islands theorem}

\noindent
 We use some elementary graph theory and apply this to
 Ahlfors five islands theorem similar to \cite{Bergweiler:2000}. If $g=(V,R)$ is a digraph, 
 $d^+(g,v)$ denotes the number of successor of same knot $v\in V$.
 The following lemma shows that there exists  a loop of length $2$
 in a digraph if each knot has sufficient many successors.

 \vspace*{12pt}
 \noindent 
 {\bf Lemma~3.1.} {\it Let $k\in \N$ and let $g=(V,R)$ be a digraph with $N:= \vert V \vert \geq 2k +1$ 
   such that $d^+(g,v) \geq N-k$ for every $v\in V$. Then there exists a $v\in V$
   and pairwise distinct knots $v_1, \ldots, v_{N-2k} \in V$ with }
   \[
    (v,v_j), (v_j,v) \in R \quad \hbox{\it for $j=1, \ldots, N-2k$.}
   \]    
 \vspace*{2pt}
  \noindent    
 {\bf Proof.} For every $v\in V$ choose pairwise distinct knots $v_1, \ldots, v_{N-k} \in V$
   with $(v,v_j) \in R$ for $j=1, \ldots, N-k$. It  follows
   \[
      \vert R \vert \geq N \cdot (N-k) = N^2 - N\cdot k. 
   \]
    Assume the claim is false. Let  $v\in V$. 
    There exists at most $N-2 k -1$ pairwise distinct knots
    $w\in \{v_1, \ldots, v_{N-k}\}$ with $(w,v) \in R$. Hence, for every $v\in V$
    there are at least $(N-k) - (N-2k -1) = k +1$ pairwise distinct knots 
    $w\in V$ with $(v,w) \in R$ and $(w,v) \not\in R$. 
    Thus $\vert R \vert \leq N^2 - N (k+1) = N^2 - Nk - N$,
    a contradiction. \, \qed

 \vspace*{12pt}
 \noindent 
 {\bf Definition~3.2.}
   Let $D$ be a domain in $\cd$, $W\subseteq \cd$ and let
   $f:D \to \cd$ be a meromorphic function. If there exists a domain $U\subseteq D$ such
   that $f\vert_{U}: U \to W$ is bijective, $U$ is called a 
   {\bf (simple) island over $W$ in $D$} (with respect to $f$).

 \vspace*{12pt}
\noindent 
   We formulate one version of the Ahlfors five islands theorem, which
   is the key tool here and one of the main tools in the complex dynamics in general.
   
 \vspace*{12pt}
   \noindent 
 {\bf Theorem~3.3.} (Ahlfors  five islands theorem \cite{Bergweiler:2000})
   Let $D_1, \ldots, D_5$ be Jordan domains on the Riemann sphere $\cd$ with pairwise
   disjoint closures and  let $D\subseteq \cd$ be a domain. 
   Then, the family of all meromorphic functions $f:D\to \cd$ with the property that 
   $f$ has no simple island over $D_j$ in $D$ for $j=1, \ldots, 5$ is normal.

 \vspace*{12pt}
   \noindent 
 {\bf Notation~} 
   The open ball centered at a point $x$ and of radius $r> 0$ will be 
   denote by $B(x,r)$.
   Let $D$ be a domain in $\cd$ and let $f:D \to \cd$ be a meromorphic function.
    For all $x = (x_1, \ldots, x_n) \in D^n$, $n\in \N$, $\gamma, \delta \in \R$ with $0<\gamma < \delta$
    and  $V_n := \{1, \ldots, n\}$, we  write 
     \[
       R_f( x,\gamma,\delta) := \{ (i,j) \in V_n \times V_n : \hbox{it exists a island over $B(x_j,\delta)$ in $B(x_i,\gamma)$} \}.
     \]
    The corresponding digraph is  given by 
    \[
       g = g_n(x,\gamma,\delta,f)  := (V_n, R_f(x,\gamma,\delta)).
    \]

 \vspace*{2pt}
 \noindent 
  The following theorem is a simple consequence of the Ahlfors five islands theorem.

 \vspace*{12pt}
 \noindent 
 {\bf Theorem~3.4.} {\it Let $D$ be a domain in $\C$, $a\in D$ and let $n\in \N_{\geq 5}$, $x_1, \ldots, x_n \in \C$,
 $\delta > 0$ with $\overline{B}(x_i,\delta) \cap \overline{B}(x_j,\delta) = \emptyset$
 for all $i\not=j$ and let $(f_k)$ be a sequence of meromorphic functions on $D$ 
 such that no subsequence is normal in $a$. Then there exists a $k\in  \N$
 and pairwise distinct $i_1, \ldots, i_{n-4} \in \{1, \ldots, n\}$ such
 that there is a island (with respect to $f_k$) over $B(x_{i_j},\delta)$ in $D$ 
 for $j=1, \ldots, n-4$. }

 \vspace*{12pt}
 \noindent 
 {\bf Proposition~3.5.} 
 Let $D$ be a domain in $\cd$, $n\in \N_{\geq 5}, x=(x_1, \ldots, x_n) \in D^n$ 
 and let $\gamma, \delta \in \R$ with $0<\gamma \leq \delta$ and
 $\overline{B}(x_j,\gamma) \subseteq D$ for $j=1, \ldots, n$ such that
 $\overline{B}(x_j,\delta) \cap \overline{B}(x_i,\delta) = \emptyset$ for all $i\not= j$.
 Let $(f_k)$ a sequence of meromorphic functions on $D$ such that every subsequence is not normal
 in $x_1, \ldots, x_n$.
 Then there exists a $k\in \N$ with
 \[
    d^+(v, g_n(x,\gamma,\delta,f_k)) \geq n -4 \quad \hbox{for all $v\in V_n$.}
 \]

 \vspace*{2pt}
  \noindent    
 {\bf Proof.} Define  $F_0 := \{f_k:k \in \N\}$ and inductively
   for $m =1, \ldots, n$:
   $F_m := \bigl\{f \in F_{m-1}: d^+\bigl( j,g_n(x,\gamma,\delta,f)\bigr) \geq n-4\bigr\}$.   
   Theorem 3.4 implies that $F_m$ is not normal in $x_1, \ldots, x_m$ for $m=1, \ldots, n$, thus
   $F_m \not= \emptyset$. If we choose  $k\in \N$ with $f_k \in F_n$,
   the proof is complete.  \, \qed 

 \section{Koebe distortion theorem and iterated function system}
 \noindent 
 {\bf Notation} (Koebe distortion theorem \cite{Schiff:93}) If $G$ is a domain in $\C$ and $K\subseteq G$ compact then
 there exists a $M=M(G,K) > 0$ such that for every univalent holomorphic 
 map $f:G\to \C$ and for all $z,w \in K$
 \[
    \frac{1}{M} \leq \frac{ \vert f'(z)\vert}{\vert f'(w) \vert} \leq M.
 \]

 \vspace*{2pt}
 \noindent 
 {\bf Proposition~4.1.} {\it Let $a\in \C, \delta > 0, U:= B(a,\delta)$ and $D := B(a, \delta/2)$.
  Set $M:= M(U, \overline{D})$ and $\gamma := \frac{\delta}{8 \cdot M}$. 
  Let  $n \in \N$ and let $D_1, \ldots, D_n$ be open subsets of $B(a,\gamma)$ with 
  $\overline{D_i} \cap\overline{D_j} = \emptyset$ for $i\not=j$.
  Let $g:U \to \cd$ be a meromorphic map such that
  \[
     g\vert_{D_j}: D_j \to U 
  \]
  is bijective for $j=1, \ldots, n$.
  Then there exists a non empty, compact subset
  $A$ of $U$ and a homeomorphism $\Phi: \Sigma_n \to A$ with $\Phi \circ \sigma_n = g \circ \Phi$,
  where $\sigma_n:\Sigma_n \to \Sigma_n$ denotes the (one-side) shift 
  on $\Sigma_n = \{1, \ldots, n\}^\N$.}

 \vspace*{12pt}
  \noindent    
 {\bf Proof.}
 Let $T_j$ the branch of $g^{-1}$ which maps $D$ to $g^{-1}(D) \cap D_j$ for $j=1, \ldots, n$.
 It follows that $T_j$ can extended to a injective function on $U$. Consider
 \[
   h: \D \to \C, \quad z \mapsto \frac{ T_j(a+\frac{\delta}{2} \cdot z) - T_j(a)}{2 \cdot \gamma}
 \]
 where $\D = \{z \in \C: \vert z \vert < 1 \}$.
 $h$ is a holomorphic map with $h(0) = 0$ and $\vert T_j(a+\frac{\delta}{2} \cdot z) - T_j(a)\vert \leq \vert T_j(a+\frac{\delta}{2} \cdot z) - a\vert + \vert T_j(a) -a\vert < 2 \cdot \gamma$
 for all $z\in \D$, thus $h(\D) \subseteq \D$. It follows $\frac{\delta}{4 \cdot \gamma} \cdot \vert T_j'(a)\vert = \vert h'(0) \vert \leq 1$
 by Schwarz's lemma. Hence
 \[
    \vert T_j'(z)\vert \leq M \cdot \vert T_j'(a) \vert \leq M \cdot \frac{4 \cdot \gamma}{\delta} = M \cdot \frac{\delta}{2 \cdot M} \cdot \frac{1}{\delta} = \frac{1}{2} \quad \hbox{for all $z\in \overline{D}$,}
 \]
 i.e. $T_j\vert_{\overline{D}}$ is a contraction and $( \overline{D},d,T_1, \ldots, T_n)$ is a hyperbolic iterated function system,
 where $d$ denotes the Euclidean metric \cite{Edgar:90, Edgar:98}. Let $A$ be the corresponding attractor and let $\Phi$ be the
 address map of this iterated function system, i.e. $\Phi(\omega) = \lim_{k\to \infty} T_{\omega_1} \circ \cdots T_{\omega_k}(x)$
 independent of the point $x \in \overline{D}$. 
  $A$ is a non-empty, compact subset of $\overline{D} \subseteq U$ and $\Phi: \Sigma_n \to A$ is a 
  homeomorphism if $\Phi$ injective. 
  $g\circ T_j = \id$ and the continuity of $g$ implies for all 
  $\omega = (\omega_k) \in \Sigma_n, x \in \overline{D}$
 \begin{eqnarray*}
   \Phi \bigl(\sigma_n (\omega)\bigr) &=& \lim_{k\to \infty} (T_{\omega_2} \circ \cdots \circ T_{\omega_{k+1}} )(x) = \lim_{k\to \infty} (g \circ T_{\omega_1}) \circ (T_{\omega_2} \circ \cdots \circ T_{\omega_{k+1}})(x) \\
        &=& g \bigl( \lim_{k\to \infty} (T_{\omega_1} \circ \cdots \circ T_{\omega_{k+1}} )(x) \bigr) =  g \bigl( \Phi (\omega) \bigr).
 \end{eqnarray*}
 Let $\omega =(\omega_k), \omega' = (\omega_k') \in \Sigma_n$
 with $\omega \not= \omega'$. Choose a minimal $k\in \N$ with $\omega_k \not= \omega_k'$
 and set $T := T_{\omega_1} \circ \cdots \circ T_{\omega_{k-1}}$ if $k\geq 2$ and $T=\id$
 otherwise. Using $T_{\omega_k} : D \to D_{\omega_k}$, $T_{\omega_k'}: D \to D_{\omega_k'}$ and
 $\overline{D}_{\omega_k} \cap \overline{D}_{\omega_k'} = \emptyset$,
 as well as the injectivity of $T$, it follows
 \begin{eqnarray*}
  (T_{\omega_1} \circ \cdots \circ T_{\omega_{k+\ell}} )(\, \overline{D}\, ) \cap (T_{\omega'_1} \circ \cdots \circ T_{\omega'_{k+\ell}} )(\,\overline{D}\,) &\subseteq& T(D_{\omega_k}) \cap T(D_{\omega'_k}) \\
     &\subseteq& T(\,\overline{D}_{\omega_k}) \cap T(\,\overline{D}_{\omega'_k}) = \emptyset.
 \end{eqnarray*}
  For $x\in \overline{D}$ the points $(T_{\omega_1} \circ \cdots \circ T_{\omega_{k+\ell}})(x)$ and
  $(T_{\omega'_1} \circ \cdots \circ T_{\omega'_{k+\ell}})(x)$ lies in two disjoint
  closed sets which are independent of $\ell \in \N$. Thus $\Phi$ is injective.  \, \qed 

 \vspace*{12pt}
 \noindent
 {\bf Theorem~4.2.} {\it Let $D$ be a open subset of $\cd$, let $F$ be a family of meromorphic
  functions on $D$ and $U \subseteq D$ open with $\infty \not\in U$. 
  Let $n\in \N_{\geq 9}$ and $x_1, \ldots, x_n \in U$ pairwise distinct such
  that $F$ is not normal in $x_1, \ldots, x_n$. Then there exists a non empty,
  compact subset $A$ of $U$ and a homeomorphism $\Phi: \Sigma_n \to A$ with
  \[
    \Phi \circ \sigma_n = f^2 \circ \Phi
  \]
  for same $f\in F$.}
  
 \vspace*{12pt}
  \noindent    
 {\bf Proof.} 
   Choose same $\delta \in \, ]0, \infty[$ with $\overline{B}(x_j,\delta) \subseteq U$ for $j=1, \ldots, n$ and
   $\overline{B}(x_j, \delta) \cap \overline{B}(x_i, \delta) = \emptyset$
   for all $i\not= j$. There exists a sequence $(f_k)$ in $F$ such that each subsequence
   is not normal in $x_1, \ldots, x_n$. Let $M := M( B(x_1, \delta), \overline{B}(x_1, \frac{\delta}{2}))$
   and $\gamma := \frac{\delta}{8 \cdot M}$, $x =(x_1, \ldots, x_n)$. 
   Choose  $k\in \N$ such that  $d^+( v, g_n(x, \gamma, \delta, f_k)) \geq n-4$
   for every $v \in V_n$ by Proposition 3.5.
   Using Lemma 3.1 there exists a $v \in V_n$ and pairwise distinct $v_1, \ldots, v_{n-8} \in V_n$ 
   with
   \[
      (v, v_j), (v_j, v) \in R := R_{f_k}(x, \gamma, \delta) \quad \hbox{for $j=1, \ldots, n-8$.}
   \]
   $(v,v_j) \in R$ implies that there exists an island $D_j$  over $B(x_{v_j},\delta)$ in $B(x_{v},\gamma)$ (with respect to $f_k$)
   for $j=1, \ldots, n-8$, and  $(v_j, v) \in R$ implies that there exists
   a island $E_j$ over $W := B(x_v,\delta)$ in $B(x_{v_j},\gamma)$ (with respect to $f_k$)
   for $j=1, \ldots, n-8$. Thus
   \[
     f_k\vert_{D_j} : D_j \to B(x_{v_j}, \delta), \quad f_k\vert_{E_j}: E_j \to W, \qquad \hbox{$j=1, \ldots, n-8$}
   \]
   are bijective. Set $D_j^* := \bigl( f_k \vert_{D_j} \bigr)^{-1}(E_j) \subseteq D_j \subseteq W$ for $j=1,\ldots, n-8$ and
   \[
     g:= f_k^2\vert_{W} : W \to \cd. 
   \]
   It follows that $g\vert_{D_j^*} : D_j^* \to W = B(x_v, \delta)$ is bijective. 
   Hence $g = f_k^2\vert_W$ is locally conjugated to the shift $\sigma_n$ on $\Sigma_n$
   by Proposition 4.1. \, \qed 

\section{Polynomial-like mappings}
\noindent
  The entropy of a rational function $Q:\cd \to \cd$ is given by 
  \cite{Ljubich:83} as $\log \deg Q$. We use this and the concept
  of polynomial-like mappings \cite{Hubbard:85} to compute 
  the entropy of entire transcendental functions.

 \vspace*{12pt}
 \noindent
 {\bf Definition~5.1} Let $U,V$ be simple connected bounded domains in $\C$
 with $\overline{U}\subseteq V$ and let $f:U \to V$. If $f$ is a proper map 
 of degree $d\in \N$, we call $(f,U,V)$ polynomial-like map (of degree $d$).

 \vspace*{12pt}
 \noindent
  The next Proposition gives a variety of examples of
  polynomial-like mappings. 

 \vspace*{12pt}
 \noindent
 {\bf Proposition~5.2} 
 Let $D$ be  a bounded domain in $\C$, $f:D\to \C$ be a holomorphic map and
 let $V$ be a open and bounded subset of $\C$. 
 Then for every connected component $U$
 of $f^{-1}(V)$ with $\overline{U} \subseteq D$ the map 
 \[
   f\vert_U:U\to V 
 \]
 is proper, i.e. a polynomial-like map if  $V$ is simple connected.
 
  \vspace*{12pt}
  \noindent    
 {\bf Proof.}
   Let $U$ be a connected component of $f^{-1}(V)$ with
   $\overline{U} \subseteq D$ and let $(z_n)$ be a sequcne
   in $U$ with $\dist(\partial U, z_n) = \inf_{z\in \partial U} \vert z - z_n\vert \to 0$.
   $\partial U$ is compact, so w.l.o.g. we assume the 
   existence of same $z \in \partial U \subseteq D$ with $z_n \to z$.
   The compactness of $\overline{V}$ implies w.l.o.g. $f(z_n) \to w$
   for same $w$, and $f(z) = w$.
   Assume $w \in V$. Let $W$ be a component of $f^{-1}(V)$ with
   $z \in W$. Thus $z \in W \cap \partial U$, i.e.
   $U\cap W\not=\emptyset$, hence $W=U$.
   It follows $z \in W=U$, a contradiction.
   So we have $w\in \partial V$. Using a characterisation
   of proper mappings \cite{Conway:96}, the proof is done. \, \qed 

 \vspace*{12pt}
 \noindent
 {\bf Theorem~5.3} (Straightening Theorem \cite{Hubbard:85}) Let $d\in \N$ and let $(f,U,V)$ be a polynomial-like map
  of degree $d$. Then there exists a polynomial $P:\C\to \C$ and a
  quasi-conform map $\varphi:\C\to \C$ with 
  \[
   \varphi \bigl( f(z) \bigr) = P \bigl( \varphi(z) \bigr) \quad \hbox{for all $z\in U$}
 \]
  such that $\varphi(U)$ contains the filled Julia set $K(P) := \{z \in \C: P^n(z) \;\not\!\!\longrightarrow \infty \}$ of $P$.
 
 \vspace*{12pt}
 \noindent
 {\bf Corollary~5.4} 
   Let $(f,U,V)$ be a polynomial-like map of degree $d\in \N$ and let
   $\varphi:\C \to \C$ with $\varphi\bigl( f(z) \bigr) = P\bigl( \varphi(z) \bigr)$
   for all $z\in U$ by Theorem 5.3. Then $f\vert_A : A \to A$,
   where $A=\varphi^{-1}\bigl( K(P) \bigr)$, and
   \[
      h(f\vert_A) = \log d.
   \]    
  \vspace*{2pt}
  \noindent    
 {\bf Proof.}
   We have $P\vert_{K(P)}: K(P) \to K(P)$
    and $A = \varphi^{-1}(K(P)) \subseteq \varphi^{-1} (\varphi(U)) = U$.
    Thus $f(z) = \varphi^{-1}\bigl( P(\varphi(z) \bigr) \in \varphi^{-1}(K(P)) = A$
    for all $z\in A$, i.e. $f\vert_A: A \to A$. 
    Proposition 1.2, Theorem 2.4  and \cite{Ljubich:83} implies
    \[
      h(f\vert_A) = h\bigl( P, K(P) \bigr) = h(P) = \log d.  \qquad   \, \qed 
    \]

\section{The main result} 
 \vspace*{12pt}
 \noindent
 {\bf Theorem~6.1.} {\it Let $f$ be a entire transcendental function. 
    Then $h(f) = \infty$.}

 \vspace*{12pt}
  \noindent    
 {\bf Proof.} Let $f_n(z) := \frac{f(n\cdot z)}{n}$ for all $z \in \C, n\in \N$.
  Then $f$ and $f_n$ are topological conjugated, thus $h(f) = h(f_n)$ for all $n\in \N$.
  Let $F:=\{ f_n : n\in \N\}$. See \cite{Schiff:93} for the definition of quasi-normality.
  \begin{enumerate}
    \item $F$ is not quasi-normal.
    
          Then we find a sequence $(f_{n_k})$ in $F$ and a sequence $(x_j)$
	  of pairwise distinct points in $\C$ such that $(f_{n_k})$ is 
	  not normal in $x_j$ for all $j\in \N$. Let $m\in \N_{\geq 9}$.
	  Theorem 4.2 implies that a non empty, compact set $A\subseteq \C$ and a
	  homeomorphism $\Phi: \Sigma_m \to A$ as well as a $k\in \N$ exists
	  with $\Phi \circ \sigma_m = f^2_{n_k} \circ \Phi$. It follows
	  $f_{n_k}^2(A) = \Phi\bigl( \sigma_m\bigl( \Phi^{-1}(A) \bigr) \bigr) = A$ and
	  \[
	     2 \cdot h(f) = 2 \cdot h(f_{n_k}) = h(f^2_{n_k}) \geq h(f_{n_k}^2, A) = h(f^2_{n_k}\vert_A ) = h(\sigma_m) = \log m,
	  \]
	    thus $h(f) = \infty$.
    
    \item $F$ is quasi-normal.
    
          Then there exists a finite set $E\subseteq \C$ and a sequence
	  $(f_{n_k})$ in $F$ such that $(f_{n_k})$ convergence local uniformly
	  on $\C \setminus E$. $\{ f_{n_k}: k \in \N\}$ is not normal in $0$,
	  thus $\vert E \vert \geq 1$ and it follows $f_{n_k} \to \infty$ on
	  $\C\setminus E$ \cite[Prop. A.2]{Schiff:93}. Choose same $r \in \, ]0, 1[$ with
	  $f_{n_k}(z) \to \infty$ for all $z\in B(0,r)$. $f_{n_k}(0) = \frac{f(0)}{n_k} \to 0$
	  implies the existence of same $k_0 \in \N$ with $\vert f_{n_k}(0)\vert < r$ and
	  \[
	     \vert f_{n_k}(z) \vert > 1 \quad \hbox{for all $k \in \N_{\geq k_0}, z \in \partial B(0,r)$.}
	  \]
          Let $k\in \N_{\geq k_0}$ and let $U_k$ be the connected component
	  of $f_{n_k}^{-1}\bigl( B(0,r) \bigr)$ which contains $0$. It follows
	  $\overline{U_k} \subseteq B(0,r)$  by $\partial B(0,r) \cap f_{n_k}^{-1}\bigl( B(0,r) \bigr) = \emptyset$.
	  Thus $\bigl(f_{n_k}, U_k, B(0,r)\bigr)$ is a polynomial-like map by Proposition 5.2.
	  Corollary 5.4 implies that it is enough to show that the sequence 
	  of the degrees of $f_{n_k}\vert_{U_k}: U_k \to B(0,r)$ is unbounded.
	  We can assume that the number of zeros of $f$ is unbounded (otherwise we
	  conjugate $f$). Let $m\in \N$ and choose $R> 0$ such that
	  $B(0,R)$ contains at least $m$ zeros of $f$. Set
	  \[
	    M:= \max_{\vert z \vert = R} \vert f(z)\vert = \max_{\vert z \vert \leq R } \vert f(z)\vert
	  \]
	  and let $W$ be the connected component of $f^{-1} \bigl( B(0,R) \bigr)$
	  which contains $0$.  
	  $B(0,R)$ is a connected set with $B(0,R) \subseteq f^{-1}\bigl( B(0,M) \bigr)$
	  which contains $0$, thus $B(0,R) \subseteq W$, i.e. $W$ contains at least $m$ zeros
	  of $f$. Choose  $k\in \N_{\geq k_0}$ with $\frac{M}{n_k} < 1$ and set
	  $g(z) := n_k \cdot z$ for all $z\in \C$. If $z\in g^{-1}\bigl( f^{-1}\bigl( B(0,M) \bigr)\bigr)$,
	  then $\vert f(n_k \cdot z)\vert < M$, thus $\vert f_{n_k}(z) \vert = \vert \frac{ f( n_k \cdot z)}{n_k} \vert < \frac{M}{n_k} < 1$,
	  hence $z\in f_{n_k}^{-1}(\D)$. It follows $g^{-1}\bigl( f^{-1}\bigl( B(0,M) \bigr) \bigr) \subseteq f_{n_k}^{-1}(\D)$,
	  especially $g^{-1}(W) \subseteq f_{n_k}^{-1}(\D)$. Thus $g^{-1}(W)$ is the
	  connected component of $g^{-1}\bigl( f^{-1}\bigl( B(0,M) \bigr) \bigr) \subseteq f_{n_k}^{-1}(\D)$
	  which contains $0$ and the maximality of $U_k$ implies $g^{-1}(W) \subseteq U_k$.
	  Let $z\in W$ with $f(z) = 0$. Then $f_{n_k}\bigl( g^{-1}(z) \bigr) = \frac{f(n_k \cdot \frac{1}{n_k} z) }{n_k} = \frac{f(z)}{n_k} = 0$
	  thus $g^{-1}(W) \subseteq U_k$ contains at least $m$ zeros of $f_{n_k}$.   \, \qed 
  \end{enumerate}

 \vspace*{12pt}
 \noindent
 {\bf Remarks~6.2} 
   The proof of Theorem 6.1 shows, that only the following 
   properties of the topological entropy where needed:
   \begin{enumerate}
     \item[(i)] The topological entropy is a topological property,
     \item[(ii)] $h(f\vert A) \leq h(f)$, where $A\subseteq X$ closed with $f(A) \subseteq A$,
     \item[(iii)] $h(f^2) \leq 2 \cdot h(f)$.
   \end{enumerate}
  So, one can replace this definition of topological entropy (on non-compact spaces)
  by an arbitrary which satisfies the obove properties.

\nonumsection{Acknowledgments}
\noindent
 The author would like to thanks Walter Bergweiler for stimulating discussions and
 contributing ideas. Also
 thanks to Klaus Schmidt and Peter Walters for their nice 
 hospitality at the university of Warwick in Summer 2003. This paper is
 based on the author's diploma thesis \cite{Wendt:2002}.

\nonumsection{References}

\end{document}